\documentclass[12pt,reqno]{amsart}
\usepackage{amssymb}
\usepackage{amsmath}
\usepackage{amsthm}
\usepackage{eucal}
\usepackage{color}
\usepackage{amsfonts}
\usepackage[pdftex]{graphicx}
\usepackage[T1]{fontenc}

\newcommand{\R}{\mathbb R}
\newcommand{\So}{\mathbb S}

\newcommand{\p}{\partial}

\newcommand{\sgn}{\text{sgn}}

\renewcommand{\a}{\alpha}
\renewcommand{\b}{\beta}

\newtheorem{theorem}{Theorem}[section]

\newtheorem{remark}[theorem]{Remark}

\numberwithin{equation}{section}
\newcommand{\be}{\begin{equation}}
\newcommand{\ee}{\end{equation}}
\newcommand{\bp}{\begin{proof}}
\newcommand{\ep}{\end{proof}}
\newcommand{\bel}{\begin{equation}\label}
\newcommand{\eeq}{\end{equation}}
\newtheorem{thm}{Theorem}[section]

\theoremstyle{remark}

\numberwithin{equation}{section}

\begin{document}
\title[CH equation]{On special properties of solutions to Camassa-Holm  equation and related models.}
\author{C. Hong}
\address[C. Hong]{Department  of Mathematics\\
University of California\\
Santa Barbara, CA 93106\\
USA.}
\email{christianhong@ucsb.edu}

\author{F. Linares}
\address[F. Linares]{Instituto de Matem\'atica Pura e Aplicada-IMPA\\
Estrada Dona Castorina, 110
Jardim Botânico | CEP 22460-320\\
Rio de Janeiro, RJ - Brasil.}
\email{linares@impa.br}

\author{G. Ponce}
\address[G. Ponce]{Department  of Mathematics\\
University of California\\
Santa Barbara, CA 93106\\
USA.}
\email{ponce@math.ucsb.edu}

\keywords{Nonlinear dispersive equation,  unique continuation, regularity}
\subjclass{Primary: 35Q35, 35B65 Secondary: 76B15}
\dedicatory{Dedicated to Professor Tohru Ozawa}

\begin{abstract} 
We study unique continuation properties of solutions to the b-family of equations. This includes the Camassa-Holm and the Degasperi-Procesi models.  We prove  that for both, the initial value problem and the periodic boundary value problem, the unique continuation results found in \cite{LiPo} are optimal. More precisely, the result established there for the constant $c_0=0$ fails for any constant $c_0\neq 0$.
\end{abstract}

\thanks
{F.L. was partially supported by CNPq grant 310329/2023-0 and FAPERJ grant E-26/200.465/2023.}

\maketitle

\section{Introduction and main results}


The  Camassa-Holm (CH) equation 
\begin{equation}\label{CH1}
\partial_tu-\partial_t \partial_x^2 u+3 u \partial_xu=2\partial_xu \partial_x^2u+ u\partial_x^3u, \hskip4pt \;t\in\R,\;x\in\R \:(\text{or}\;\So),
\end{equation}
was first noted by Fuchssteiner and Fokas \cite{FF} in their work on hereditary symmetries. Later, it was derived  as a model for shallow water waves  by Camassa and Holm \cite{CH}, who also examined its solutions. It has also appeared as a model in nonlinear dispersive waves in hyper-elastic rods, see \cite{Dai}, \cite{DaHu}. 

The CH equation \eqref{CH1} has received extensive attention due to its remarkable properties, among them the fact that it is a bi-Hamiltonian completely integrable model (see \cite{BSS},  \cite{CH}, \cite{CoMc}, \cite{Mc}, \cite{M}, \cite{Par} and references therein).

 The CH equation is a member of the so called b-family derived in \cite{HoSt}
 \begin{equation}\label{B-fam}
\partial_tu-\partial_t \partial_x^2 u+(b+1) u \partial_xu=b\partial_xu \partial_x^2u+ u\partial_x^3u, \hskip10pt b\in\R.
\end{equation}
 
This family of equations can be written as 
\begin{equation}\label{b-fam}
\partial_tu+u\partial_x u +\,\partial_x(1-\partial_x^2)^{-1}\Big(\,\frac{b}{2}u^2  +\frac{3-b}{2} (\partial_xu)^2\Big)=0.
\end{equation}

For $b=2$ one gets the CH equation and for  $b=3$ the Degasperi-Procesi (DP) equation \cite{DP}, the only bi-hamiltonian and integrable models in this family, see \cite{Iv}, \cite{Ma}.

 The b-family  possess ``peakon'' solutions. The  single peakon in $\R$ is explicitly given by the formula 
\begin{equation}\label{peakon}
u_c(x,t)=c\, e^{-|x-ct|}, \hskip20pt c\in \R.
\end{equation}

The initial value problem (IVP)  as well as the initial periodic boundary value problem (IPBVP)  associated to the equation \eqref{b-fam} has been extensively examined. 
In particular, in \cite{LiOl} and \cite{Ro}  strong local well-posedness (LWP) of the IVP associated to the CH equation was established in the Sobolev space
\[
H^s(\R)=(1-\partial_x^2)^{-s/2}L^2(\R),\;\;\;\;\;\;s>3/2.
\]
The argument in \cite{LiOl} and \cite{Ro} extends, without a major modification, to all the equations in the b-family, and to the IPBVP associated to them.

\begin{thm}[\cite{LiOl}, \cite{Ro}]\label{thm0-ch}\hskip15pt

\begin{enumerate}
\item
Let   $s>3/2$. For any $u_0\!\in\!  H^s(\R)$, there exist  $T\!=\!T(\|u_0\|_{s,2})>0$ and a unique solution $u=u(x,t)$ of the IVP associated to the b-family of equations \eqref{b-fam}
such that
\begin{equation}
\label{class-sol}
u\in
C([-T,T]\!:\!H^s(\R))\cap C^1((-T,T)\!:\!H^{s-1}(\R)).
\end{equation}
Moreover, the map $u_0\mapsto u$, taking the data to the solution, is locally continuous from $H^{s}(\R)$ into $C([-T,T]\!:\!H^s(\R))$.
\vskip.1in
\item
Let   $s>3/2$. For any $u_0\!\in\!  H^s(\So)$, there exist  $T\!=\!T(\|u_0\|_{s,2})>0$ and a unique solution $u=u(x,t)$ of the IPBVP associated to the b-family of equations \eqref{b-fam}
such that
\begin{equation}
\label{class-sol1}
u\in
C([-T,T]\!:\!H^s(\So))\cap C^1((-T,T)\!:\!H^{s-1}(\So)).
\end{equation}
Moreover, the map $u_0\mapsto u$, taking the data to the solution, is locally continuous from $H^{s}(\So)$ into $C([-T,T]\!:\!H^s(\So))$.
\end{enumerate}

\end{thm}


One observes that the peakon solutions do not belong to these spaces. In fact, 
\begin{equation}
\label{space}
\phi(x)=e^{-|x|}\notin W^{p,1+1/p}(\R)\;\;\;\;\;\text{for any}\;\;\;\;\;p\in[1,\infty),
\end{equation}
where $W^{s,p}(\R)=(1-\partial_x^2)^{-s/2}L^p(\R)$ with $\,W^{s,2}(\R)=H^s(\R)$. However,
\[
\phi(x)=e^{-|x|}\in W^{1,\infty}(\R),
\]
where $W^{1,\infty}(\R)$ denotes the space of Lipschitz functions. In this regards, the following weaker version of  the LWP for  the IPBVP associated to the CH equation was given in \cite{LKT}. This was extended in \cite{LiPoSi} to the case of the IVP. The proofs there apply to the b-family \eqref{b-fam}.

\begin{thm}[\cite{LKT}, \cite{LiPoSi}]\label{thm1-ch}\hskip15pt

\begin{enumerate}
\item
Given  $u_0\in X\!\equiv\! H^1(\R)\cap W^{1,\infty}(\R)$, there exist  $T\!=\!T(\|u_0\|_X)>0$ and a unique solution $u=u(x,t)$ of the IVP associated to the b-family equations \eqref{b-fam}
such that
\begin{equation}
\begin{aligned}
\label{class-sola15}
u\in
&C([-T,T]\!:\!H^1(\R))\cap C^1((-T,T)\!:\!L^2(\R))\\
&\hskip5pt\cap L^{\infty}([-T,T]\!:\!W^{1,\infty}(\R))\equiv Z_T\cap L^{\infty}([-T,T]\!:\!W^{1,\infty}(\R)).
\end{aligned}
\end{equation}
Moreover, the map $u_0\mapsto u$, taking the data to the solution, is locally continuous from $X$
into $Z_{T}$.

\vskip.1in

\item
Given  $u_0\in X\!\equiv\! H^1(\So)\cap W^{1,\infty}(\So)$, there exist  $T=T(\|u_0\|_X)>0$ and a unique solution $u=u(x,t)$ of the IPBVP associated to the b-family equations \eqref{b-fam}
such that
\begin{equation}
\begin{aligned}
\label{class-sola151}
u\in
&C([-T,T]\!:\!H^1(\So))\cap C^1((-T,T)\!:\!L^2(\So))\\
&\hskip5pt \cap L^{\infty}([-T,T]\!:\!W^{1,\infty}(\So))\equiv Z_T\cap L^{\infty}([-T,T]\!:\!W^{1,\infty}(\So)).
\end{aligned}
\end{equation}
Moreover, the map $u_0\mapsto u$, taking the data to the solution, is locally continuous from $X$
into $Z_{T}$.

\end{enumerate}
\end{thm}

For further results regarding the existence and uniqueness of solutions to equations in  the b-family we refer to \cite{BC, BCZ, CocKa, CoEs1, CoEs2, CoEs3, CoMo, EY, GHR1, GHR2, HHG, HKM, M} and references therein.

\medskip

Our main interest here is concerned with unique continuation properties for solutions of the b-family \eqref{b-fam}.  We shall consider  both, the IVP and the IPBVP  associated to the b-family \eqref{b-fam}. In this regards, we have a result proved in \cite{LiPo} and slightly improved in \cite{HP}:

\begin{thm} [\cite{LiPo}]\label{IVPCH2}\hskip15pt

\begin{enumerate}
\item
Let $\,u=u(x,t)$ be a solution of the IVP associated to the b-family of equations \eqref{b-fam}, with $b\in(0,3]$, in the class described in part (1) of Theorem \ref{thm0-ch} or Theorem \ref{thm1-ch}.
If there exist $t_0\in(-T,T)$ and $\alpha,\beta\in\R$ with $ \alpha<\beta$, such that
\begin{equation}
\label{cond1-ch}
u(x,t_0)=0,\;\;\;\;\;\;\;x\in[\alpha,\beta],
\end{equation}
and
\be
\label{cond2-ch}
\partial_tu(x,t_0)\;\;\;\;\text{is not strictly decreasing on}\;\;\;\;(\alpha,\beta),
\end{equation}
then $\,u\equiv 0$.

\vskip.1in
\item
Let $\,u=u(x,t)$ be a solution of the IPBVP associated to the b-family of equations \eqref{b-fam}, with $b\in(0,3]$, in the class described in part (2) of Theorem \ref{thm0-ch} or Theorem \ref{thm1-ch}.
If there exist $t_0\in(-T,T)$ and $\alpha,\beta\in (0,1)\sim\So$ with $ \alpha<\beta$, such that
\begin{equation}
\label{cond1-ch1}
u(x,t_0)=0,\;\;\;\;\;\;\;x\in[\alpha,\beta],
\end{equation}
and
\be
\label{cond2-ch1}
\partial_tu(x,t_0)\;\;\;\;\text{is not strictly decreasing on}\;\;\;\;(\alpha,\beta),
\end{equation}
then $\,u\equiv 0$.
\end{enumerate}

\end{thm}

\begin{remark}\hskip15pt

\begin{enumerate}
\item
From hypothesis \eqref{cond1-ch} (or \eqref{cond1-ch1}) and the equation, it follows that $\partial_tu(\cdot,t_0)$ is continuous in $[\alpha,\beta]$. Hence the point-wise evaluation in \eqref{cond2-ch} (or \eqref{cond2-ch1}) makes sense. 
\item
In the case of the Korteweg-de Vries (KdV)  equation, see \cite{KdV},
$$
\partial_tu+\partial_x^3u+u\partial_xu=0,
$$
(or any other dispersive local model) part (1) Theorem \ref{IVPCH2} fails. In this case, for a regular enough solution $u=u(x,t)$ of the KdV,  the hypothesis \eqref{cond1-ch} implies that in \eqref{cond2-ch}. Thus, in the case of the KdV, one needs a stronger hypothesis. For example, in \cite{SaSc} it was shown that the conclusion of Theorem \ref{IVPCH2} holds  if one assumes \eqref{cond1-ch}  on an open set of $\,\R\times[-T,T]$ (or of $\So\times[-T,T]$. However, this result in \cite{SaSc} applies to any two solutions 
$u_1,\,u_2$ of the KdV, but here, for the CH equation, we are assuming $u_2\equiv 0$.
\item 
The hypothesis \eqref{cond1-ch} and \eqref{cond2-ch}, restricted just to an interval, also appear in \cite{KPV19} for the case of the Benjamin-Ono equation (\cite{Be}, \cite{On})and in \cite{HP} for the Benjamin-Bona-Mahony equation (\cite{BBM}), both non-local models.
\item 
In the case of the IPBVP, a different kind of unique continuation results were established in \cite{BrCo}.

\end{enumerate}
\end{remark}


We shall prove in both cases, for the IVP and the IPBVP, that the hypotheses \eqref{cond1-ch} and \eqref{cond1-ch1} in Theorem \ref{IVPCH2} are optimal. More precisely, the results fail for any  $c_0\in \R-\{0\}$.

\begin{thm} \label{IVPCH3}

Let $b\in(0,3]$. Given any $c_0\in \R-\{0\}$ any $\alpha, \beta\in \R$ with $\alpha<\beta$,  there exists $u_0\in C^{\infty}_0(\R)$ with 
\begin{equation}
\label{cond3-ch}
u_0(x)=c_0,\;\;\;\;\;\;\;x\in[\alpha,\beta],
\end{equation}
such that the corresponding  solution $u=u(x,t)$ of the IVP associated to the b-equation \eqref{b-fam}
satisfies
\be
\label{cond4-ch}
\partial_tu(x,0)=0\;\;\;\;\;\;\;x\in[\alpha,\beta].
\end{equation}
In particular,   $\,u(\cdot,t)\not\equiv c_0$ for all $t\in (-T,T)$.
\end{thm}

\vskip.1in
\begin{thm} \label{IVPCH4}
Let $b\in(0,3]$. Given any $c_0\in \R-\{0\}$ and any $\alpha, \beta\in (0,1)\sim\So$ with $\alpha<\beta$,  there exists $u_0\in C^{\infty}(\So)$ with $u_0\not\equiv c_0$ and
\begin{equation}
\label{cond3b-ch}
u_0(x)=c_0,\;\;\;\;\;\;\;x\in[\alpha,\beta],
\end{equation}
such that the corresponding  solution $u=u(x,t)$ of the IPBVP associated to the b-equation \eqref{b-fam}
satisfies
\be
\label{cond4b-ch}
\partial_tu(x,0)=0\;\;\;\;\;\;\;x\in[\alpha,\beta].
\end{equation}
In particular,    $\,u(\cdot,t)\not\equiv c_0$ for all $t\in (-T,T)$.

\end{thm}

\begin{remark}\hskip15pt

\begin{enumerate}
\item

  From the hypothesis \eqref{cond3-ch} (or \eqref{cond3b-ch}) one needs to find and extension of $u_0$ to $\R$ (or to $\So$) different to $u_0\equiv c_0$ such that \eqref{cond4-ch} (or \eqref{cond4b-ch}) holds. From the equation, this is possible due to the explicit  form of the operator $\partial_x(1-\partial_x^2)^{-1}$ in $\R$ and in $\So$.
 
 \item
 If one only considers solutions in the class obtained in Theorem \ref{thm1-ch}  a simpler proof can be given which applies to both, the IVP and the IPBVP. However, if one requires higher regularity, one needs to consider each case separately. Thus, for the case of the IVP, Theorem \ref{IVPCH3}, we shall use an explicit constructive argument. The proof of Theorem \ref{IVPCH4}  is based on an existence argument which uses the implicit function theorem.

\item 

In particular, Theorem \ref{IVPCH3} implies that the \it stumpons \rm solutions, described in \cite{Le} for the CH equation, do not belong to the class  described in Theorem \ref{thm1-ch}.

\end{enumerate}
\end{remark}

\medskip

This paper is organized as follows: In Section 2 we will prove Theorem  \ref{IVPCH3}.  The proof of Theorem \ref{IVPCH4} will be
given in Section 3.

\section{Proof of Theorem \ref{IVPCH3}}



\begin{proof} [Proof of Theorem \ref{IVPCH3}]

Let us consider the case $b\in (0,3)$. The proof for the  $b=3$ will be given latter.

Using  the equation \eqref{b-fam} we need to find an extension of $u_0$, defined in the interval $x\in [\alpha,\beta]$ as $u_0(x)=c_0$, to the whole real line such that $u_0\in C_0^{\infty}(\R)$ and 
\begin{equation}
\label{a1}
\partial_tu(x,0)=-\partial_x(1-\partial_x^2)^{-1}F_b(u(x,0))=0,\;\;\;\;x\in[\alpha,\beta],
\end{equation}
where
\begin{equation}
\label{a2}
F_b(w(x,t))\equiv \Big(\,\frac{b}{2}w^2  +\frac{3-b}{2} (\partial_xw)^2\Big)(x,t),
\end{equation}
and
\begin{equation}
\label{a2a}
\begin{aligned}
\partial_x(1-\partial_x^2)^{-1}g(x)&=-\frac{1}{2}\sgn(\cdot)e^{-|\cdot|}\ast g(x)\\
&=-\frac{1}{2}\int_{-\infty}^{\infty}\sgn(x-y)e^{-|x-y|}g(y)dy.
\end{aligned}
\end{equation}
Thus, for $x\in [\alpha,\beta]$ one has that 
\begin{equation}
\label{a3}
\begin{aligned}
&-\!2\partial_x(1-\partial_x^2)^{-1}\!F_b(u(x,0))\!=\!\int^{\infty}_{-\infty}\!\!\sgn(x-y)e^{-|x-y|} F_b(u(y,0))dy\\
&=\int_{-\infty}^xe^{-x+y} F_b(u(y,0))dy-\int_x^{\infty}e^{-y+x}F_b(u(y,0))dy \\
&=e^{-x}(\int_{-\infty}^{\alpha}e^yF_b(u(y,0))dy+\int_{\alpha}^xe^yF_b(u(y,0))dy)\\
&\hskip10pt -e^x(\int_x^{\beta}e^{-y}F_b(u(y,0))dy+\int_{\beta}^{\infty}e^{-y}F_b(u(y,0))dy)\\
&=e^{-x}\int_{-\infty}^{\alpha}e^yF_b(u(y,0))dy+e^{-x}\,\frac{b\,c_0^2}{2}(e^x-e^{\alpha})\\
&\hskip10pt +e^x\,\frac{b\,c_0^2}{2}(e^{-\beta}-e^{-x})-e^x\int_{\beta}^{\infty}e^{-y}F_b(u(y,0))dy\\
&=e^{-x}A_b(u_0)+\frac{b\,c_0^2}{2}(1-e^{-x}e^{\alpha})-e^xB_b(u_0)+\frac{b\,c_0^2}{2}(e^{x}e^{-\beta}-1)\\
&=e^{-x}(A_b(u_0)-e^{\alpha}\,\frac{b\,c_0^2}{2})-e^x(B_b(u_0)-e^{-\beta}\,\frac{b\,c_0^2}{2}),
\end{aligned}
\end{equation}
where
\begin{equation}
\label{a4}
A_b(w)\equiv \int_{-\infty}^{\alpha}e^yF_b(w(y))dy,\hskip20pt B_b(w)\equiv \int_{\beta}^{\infty}e^{-y}F_b(w(y))dy.
\end{equation}
Thus, from \eqref{a1}, our problem reduces to find an extension of $u_0$ to $(-\infty,\alpha)\cup(\beta,\infty)$ with $u_0\in C^{\infty}_0(\R)$ such that 
in the interval $(-\infty,\alpha)$ one has that
\begin{equation}
\label{a5}
A_b(u_0)=\int_{-\infty}^{\alpha}e^yF_b(u_0(y))dy=e^{\alpha}\;\frac{b \,c_0^2}{2},
\end{equation}
and in the interval $(\beta, \infty)$ one has that
\begin{equation}
\label{a6}
B_b(u_0)=\int_{\beta}^{\infty}e^{-y}F_b(u_0(y))dy=e^{-\beta}\;\frac{b \,c_0^2}{2}.
\end{equation}

Notice that if $u_0(y)\equiv c_0$, one has that $F_b(u_0)=b\,c_0^2/2$ and 
\begin{equation}
\label{6a}
A_b(c_0)=e^{\alpha}\;\frac{b \,c_0^2}{2},\hskip20pt  B_b(c_0)=e^{-\beta}\;\frac{b \,c_0^2}{2},
\end{equation}
 but in this case $u_0\notin C^{\infty}_0(\R)$.

Since the arguments to show \eqref{a5} and \eqref{a6} are similar, we shall restrict ourselves just to prove \eqref{a6}. Thus, we are considering the extension of $u_0$ to the interval $(\beta,\infty)$.

First, for $\eta>0$ fixed, we define
\begin{equation}
\label{a7}
w_0(y)=
\begin{cases}
\begin{aligned}
&0,\;\;&y<\alpha-\eta-\frac{\pi}{2\gamma},\\
&c_0\cos(\gamma ((\alpha-\eta)-y)),\;\;&y\in[\alpha-\eta-\frac{\pi}{2\gamma}, \alpha-\eta],\\
&c_0,\;\;\;&y\in(\alpha-\eta,\beta+\eta),\\
&c_0\cos(\gamma (y-(\beta+\eta))),\;\;&y\in[\beta+\eta, \beta+\eta+\frac{\pi}{2\gamma}],\\
&0,\;\;&y>\beta+\eta+\frac{\pi}{2\gamma},
\end{aligned}
\end{cases}
\end{equation}
with
\begin{equation}
\label{a8}
\gamma=\sqrt{\frac{b}{3-b}}.
\end{equation}

Then,  $w_0$ is Lipchitz, and has compact support. We recall that we are trying to establish \eqref{a6} which depends only on the extension of $u_0$ to the interval $[\beta,\infty)$. Also, for $x\in [\alpha,\beta+\eta+\pi/2\gamma)$
\begin{equation}
\label{a9}
F_b(w_0(x))=\frac{b}{2}(w_0(x))^2+\frac{3-b}{2}(\partial_xw_0(x))^2=\frac{b\,c^2_0}{2},
\end{equation}
and for $x\in [\beta+\eta+\pi/2\gamma,\infty)$ 
\begin{equation}
\label{a10}
F_b(w_0(x))=\frac{b}{2}(w_0(x))^2+\frac{3-b}{2}(\partial_xw_0(x))^2=0.
\end{equation}
Therefore,
\begin{equation}
\label{a11}
F_b(w_0(x))=0<\frac{b \,c_0^2}{2}=F_b(c_0),\;\;x\in [\beta+\eta+\pi/2\gamma,\infty).
\end{equation}

Hence, there exists $c^*=c^*(b, c_0,\beta, \eta)>0$ such that  
\begin{equation}
\label{a12}
\begin{aligned}
B_b(c_0)&= \int_{\beta}^{\infty}e^{-y}F_b(c_0)dy=e^{-\beta}\,\frac{b\,c_0^2}{2}\\
&=B_b(w_0(y))+c^*= \int_{\beta}^{\beta+\eta+\pi/2\gamma}e^{-y}F_b(w_0)dy+c^*,
\end{aligned}
\end{equation}
i.e. $c^*=\int_{\beta+\eta+\pi/2\gamma}^{\infty}e^{-y}F_b(c_0)dy$.

\medskip

Now, let $\varphi \in C^{\infty}_0(\R)$ even, non-negative, with $supp \,\varphi\subset (-1,1)$ and $\int \varphi(x)dx=1$.  
Define $\varphi_{\epsilon}(x)=\frac{1}{\epsilon}\varphi(\epsilon x)$ for $\epsilon>0$ and
\begin{equation}
\label{a13}
w_{\epsilon}(x)=\varphi_{\epsilon}\ast w_0(x).
\end{equation}

Thus, $w_{\epsilon}\in C_0^{\infty}(\R)$ vanishing for $x>\beta+\eta+\frac{\pi}{2\gamma}+\epsilon$ (and $x<\alpha-\eta-\frac{\pi}{2\gamma}-\epsilon$) with $w_{\epsilon}(x)=c_0$ for $x\in [\alpha,\beta]$ if $\epsilon<\eta$. Moreover,
\begin{equation}\label{a14}
\lim_{\epsilon\downarrow 0}\int_{\beta}^{\infty} ((w_{\epsilon}-w_0)^2+(\partial_xw_{\epsilon}-\partial_xw_0)^2)(x)dx=0.
\end{equation}
 
 In fact, since for $\epsilon\in[0,1]$, $w_{\epsilon}$ vanishes for $x>\beta+\eta+\frac{\pi}{2\gamma}+1$, one has for any fixed $b\in(0,3)$ that
 $$
\int_{\beta}^{\infty} ((w_{\epsilon})^2+(\partial_xw_{\epsilon})^2)(x)dx\sim B_b(w_{\epsilon}).
$$
From \eqref{a14}, the hypothesis $b\in(0,3)$, and \eqref{a12}, it follows that
\begin{equation}
\label{a15}
\lim_{\epsilon\downarrow 0} B_b(w_{\epsilon})=B_b(w_0)<e^{-\beta}\,\frac{b\,c_0^2}{2}.
\end{equation}

Hence, fixing $\epsilon_0$ sufficiently small, $\epsilon_0<\eta$, we have : $w_{\epsilon_0}\in C_0^{\infty}(\R)$ vanishing for $x>\beta+\eta+\frac{\pi}{2\gamma}+\epsilon_0$ with $w_{\epsilon_0}(x)=c_0$ for $x\in [\alpha,\beta]$ and
\begin{equation}
\label{a16}
B_b(w_{\epsilon_0})<e^{-\beta}\,\frac{b\,c_0^2}{2}.
\end{equation}

We recall that we are looking for a function $u_0$ defined in $[\alpha,\infty)$ vanishing for $x>M$ for some $M>0$ such that is equal to $c_0$ for 
$x\in [\alpha,\beta]$ and $B(u_0)=e^{-\beta}\,\frac{b\,c_0^2}{2}.$ Thus, for any $\lambda>0$ we consider
\begin{equation}
\label{a17}
v_{\lambda}(x)=w_{\epsilon_0}(x)+\lambda \varphi(x-L),\;\;\;\; L=\beta+\eta+2\pi/\gamma+\epsilon_0+2.
\end{equation}
Hence, the supports of $w_{\epsilon_0}$ and $\varphi(\cdot-L)$ are disjoint and 
\begin{equation}
\label{a15b}
\lim_{\lambda\uparrow \infty} B_b(v_{\lambda})=\infty, \;\;\;\;\text{with}\;\;\;\;\;B_b(v_{0})=B_b(w_{\epsilon_0})<e^{-\beta}\,\frac{b\,c_0^2}{2}.
\end{equation}

Defining $G(\lambda)\equiv B_b(v_{\lambda})$, by  the continuity of  $G(\lambda)$  in $\lambda\in[0,\infty)$ and \eqref{a15}, there exists $\lambda_0>0$ such that $B_b(v_{\lambda_0})=e^{-\beta}\,\frac{b\,c_0^2}{2}.$

Finally, taking for $x>\alpha$
\begin{equation}
\label{a18}
u_0(x)=v_{\lambda_0}(x)=w_{\epsilon_0}(x)+\lambda_0 \varphi(x-L),\;\;\; L=\beta+\eta+2\pi/\gamma+\epsilon_0+2,
\end{equation} we complete the proof of \eqref{a6}. As we remarked before the proof of \eqref{a5} is similar. This concludes the proof of the theorem for the case $b\in (0,3)$.

\vskip.1in

Next, we study the  case $b=3$, i.e. for  the DP equation.This  is simpler since in this case $F_3$, defined in \eqref{a2}, does not depend on the values of the derivatives of $u$. So the desired extension of 
$u_0$ to the interval $[\beta,\infty)$ can be achieved by taking
\begin{equation}
\label{a19}
u_0(x)=v_0(x)+\lambda_0 \varphi(x-L),\;\;\; L=\beta+2,
\end{equation}
where $v_0\in C^{\infty}(\R)$ with $v_0(x)=c_0$ for $x\leq \beta$, decreasing for $x\in (\beta,\beta+1)$ if $c_0>0$ or increasing for $x\in (\beta,\beta+1)$ if $c_0<0$,  vanishing for $x\geq \beta+1$, so
\begin{equation}
\label{a20}
B_3(v_0)=\int_{\beta}^{\infty}e^{-y}F_3(v_0(y))dy\leq e^{-\beta}\,\frac{3\,c_0^2}{2},
\end{equation}
with $\varphi$ is as above and $\lambda_0$ such that 
\begin{equation}
\label{a21}
B_3(u_0)=\int_{\beta}^{\infty}e^{-y}F_3(u_0(y))dy=e^{-\beta}\,\frac{3\,c_0^2}{2}.
\end{equation}

\end{proof}
\section{Proof of Theorem \ref{IVPCH4}}

\begin{proof}[Proof of Theorem \ref{IVPCH4}]

Let $\alpha, \beta \in (0,1)$ with $[0, 1) \cong \mathbb S$. Using the equation \eqref{b-fam}, we need to find an extension of $u_0$, defined in the interval $x \in [\alpha, \beta]$ as $u_0(x) = c_0$, to the all of $\mathbb S$ such that $u_0 \in C^\infty (\mathbb S)$ and
\be
\label{b1}
\p_t u(x,0)  = - \p_x(1-\p_x^2)^{-1}F_b(u(x,0)) = 0, \quad x\in[\alpha,\beta]
\ee

where

\be
\label{b2}
F_b(w(x,t))= \left(\frac b 2 w^2 + \frac{3-b} 2 (\p_x w)^2 \right),
\ee
and
\be
\label{b2b}
\begin{aligned}
\p_x(1-\p_x^2)^{-1}g(x)&=\frac{\sinh(\cdot-[\cdot]-1/2)}{2\sinh(1/2)}\ast g(x)\\
&=\int_0^1\frac{\sinh(x-y-[x-y]-1/2)}{2\sinh(1/2)} g(y) dy,
\end{aligned}
\ee
where $[\cdot]$ denotes the greatest integer function. Thus, for $x \in [\alpha,\beta]$ one has that

\be
\label{b3}
\begin{aligned}
    & 2\sinh(1/2)\p_x (1-\p_x^2)^{-1}F_b(u(x,0))\\
    & = \int_0^1 \sinh (x - y - [x - y] - 1/2) F_b(u_0) \, dy \\
    & = \int_0^{\a} \sinh (x - y - 1/2) F_b(u_0) \, dy \\
    &\hskip10pt+ \int_{\a}^x \sinh (x - y  - 1/2) \frac {c_0^2 b} 2 \, dy \\
    &\hskip10pt +\int_x^{\b} \sinh (x - y +1/2) \frac {c_0^2 b} 2 \, dy \\
    &\hskip10pt+ \int_{\b}^1 \sinh (x - y + 1/2) F_b(u_0) \, dy \\
    & = \mbox I + \mbox {II} + \mbox {III} + \mbox {IV}
\end{aligned}
\ee

Using the formulas for $\sinh (x+y)$ and $\cosh (x+y)$ and integrating, we get

\be
\label{b4}
\begin{aligned}
    \mbox{I} &= \sinh(x-1/2) \int_0^{\a} \cosh(y) F_b(u_0) \, dy\\
    &\hskip10pt - \cosh ( x - 1/2) \int_0^{\a} \sinh(y) F_b(u_0) \, dy \\
    &= \sinh(x-1/2)A_c(u_0) - \cosh(x-1/2) A_s(u_0),
\end{aligned}
\ee

\be
\label{b5}
    \mbox{II} = \frac {c_0^2 b} 2 \left( \cosh(x - \a - 1/2) - \cosh(1/2) \right),
\ee

\be
\label{b6}
    \mbox{III} = \frac {c_0^2 b} 2 \left( \cosh(1/2) - \cosh( x- \b + 1/2) \right),
\ee
and 
\be
\label{b7}
\begin{aligned}
    \mbox{IV} &= \sinh(x+1/2) \int_{\b}^1 \cosh(y) F_b(u_0) \, dy \\
    &\hskip10pt- \cosh ( x + 1/2) \int_{\b}^1 \sinh(y) F_b(u_0) \, dy \\
    &= \sinh(x+1/2)B_c(u_0) - \cosh(x+1/2) B_s(u_0).
\end{aligned}
\ee

For convenience, let $\mu = \frac {c_0^2 b} 2$. Now consider the situation when $u_0 \equiv c_0$. This will give us 
\be
\label{b8}
\begin{aligned}
    A_c(c_0) &= \mu \int_0^{\a} \cosh(y) \, dy = \mu \sinh(\a), \\
    A_s(c_0) &= \mu \int_0^{\a} \sinh(y) \, dy = \mu (\cosh(\a) - 1), \\
    B_c(c_0) &= \mu \int_{\b}^1 \cosh(y) \, dy = \mu (\sinh(1) - \sinh(\b)) ,\\
    B_s(c_0) &= \mu \int_{\b}^1 \sinh(y) \, dy = \mu (\cosh(1) - \cosh(\b))
\end{aligned}
\ee

Plugging these values back into \eqref{b4} and \eqref{b7} and applying the $\sinh x$ and $\cosh x$ addition formulas will give us
\be
\p_t u(x,0) = \mbox I + \mbox {II} + \mbox {III} + \mbox {IV} = 0,
\ee

for $x \in [\alpha, \beta]$.

Now, let $u_0 = c_0 + \phi$ where $\phi \in C_0^\infty(\R)$, and $supp \, \phi \subset (0, \alpha) \cup (\beta, 1)$. Our problem reduces to find nontrivial $\phi$ such that in the interval $(0, \alpha)$ one has that

\be
\label{b9}
\begin{aligned}
    A_c(u_0) &= \int_0^{\a} \cosh(y)F_b(c_0 + \phi) \, dy = \mu \sinh(\a), \\
    A_s(u_0) &= \int_0^{\a} \sinh(y)F_b(c_0 + \phi) \, dy = \mu (\cosh(\a) - 1),
\end{aligned}
\ee

and in the interval $(\beta, 1)$ one has that

\be
\label{b10}
\begin{aligned}
    B_c(u_0) &= \int_{\b}^1 \cosh(y)F_b(c_0 + \phi) \, dy = \mu (\sinh(1) - \sinh(\b)), \\
    B_s(u_0) &= \int_{\b}^1 \sinh(y)F_b(c_0 + \phi) \, dy = \mu (\cosh(1) - \cosh(\b)),
\end{aligned}
\ee

with $supp \, \phi \subset (0, \alpha) \cup (\beta, 1)$.

\medskip

The arguments to show \eqref{b9} and \eqref{b10} are similar, so we will just prove \eqref{b9}. Using the definitions $\cosh x = \frac{e^x + e^{-x}} 2$ and $\sinh x= \frac{e^x - e^{-x}} 2$, we can rewrite \eqref{b9} as 

\be
\label{b11}
\begin{aligned}
    C_1(c_0 + \phi) = \int_0^{\a} e^y \left[F_b(c_0 + \phi) - F_b(c_0) \right] \, dy &= 0, \\
    C_2(c_0 + \phi) = \int_0^{\a} e^{-y} \left[F_b(c_0 + \phi) - F_b(c_0) \right] \, dy &= 0,
\end{aligned}
\ee

Now, let $\phi_1 \in C_0^\infty(\R)$ be nonnegative with $supp \, \phi_1 \subset (0, \a/3)$. Then let $\phi_2(x) = \phi_1(x - \a/3)$, and $\phi_3(x) = \phi_1(x - 2\a/3)$.  Define the function $L:\R^3\to \R^2$ as
\be
\label{b12}
\begin{aligned}
    &L(\lambda_1, \lambda_2, \lambda_3) = (L_1(\lambda_1, \lambda_2, \lambda_3), L_2(\lambda_1, \lambda_2, \lambda_3))\\
    & =(C_1, C_2)(c_0 + \lambda_1\phi_1 + \lambda_2 \phi_2 + \lambda_3 \phi_3).
\end{aligned}
\ee

Notice that $\phi_1, \phi_2, \phi_3$ all have disjoint support. This fact and \eqref{b2} applied to \eqref{b12} give us
\be
\label{b13}
\begin{aligned}
    L_1  &= \int_0^{\a} e^y \frac b 2 \left( 2c_0( \sum_{j=1}^3 \lambda_j \phi_j ) +c_0^2+\sum_{j=1}^3 \lambda_j^2 \phi_j^2  \right) \\
    &\hskip15pt+ e^y\frac {3-b} 2 \left( \lambda_1^2 (\phi_1')^2 + \lambda_2^2 (\phi_2')^2 + \lambda_3^2 (\phi_3')^2\right) \, dy   \\
    L_2  &=  \int_0^{\a} e^{-y} \frac b 2 \left( 2c_0 ( \sum_{j=1}^3 \lambda_j \phi_j ) +c_0^2+\sum_{j=1}^3 \lambda_j^2 \phi_j^2  \right) \\
    &\hskip15pt+ e^{-y}\frac {3-b} 2 \left( \lambda_1^2 (\phi_1')^2 + \lambda_2^2 (\phi_2')^2 + \lambda_3^2 (\phi_3')^2\right)\, dy.
\end{aligned}
\ee

Note that finding a $(\lambda_1, \lambda_2, \lambda_3) \ne(0, 0, 0)$ that solves $L(\lambda_1, \lambda_2, \lambda_3) = (0, 0)$ gives the desired result. Our goal is to apply the Implicit Function Theorem to get a curve parametrized in $\lambda_3$  such that $L(\lambda_1(\lambda_3), \lambda_2(\lambda_3), \lambda_3) = (0, 0)$ is satisfied. This requires us to show that 

\be
\label{b14}
\det \begin{bmatrix}
    \frac {\partial L_1}{\partial \lambda_1} & \frac{\partial L_1} {\partial \lambda_2} \\
    \frac {\partial L_2}{\partial \lambda_1} & \frac{\partial L_2} {\partial \lambda_2}
\end{bmatrix}(0,0,0) \ne 0.
\ee

Computing the partial derivatives one gets
\be
\label{b15}
\begin{aligned}
    \frac {\partial L_1}{\partial \lambda_1}(0, 0, 0) &= b c_0 \int_0^{\a} e^y \phi_1(y) \, dy, \\
    \frac {\partial L_1}{\partial \lambda_2}(0, 0, 0) &=  b  c_0\int_0^{\a} e^y \phi_2(y) \, dy ,\\
    \frac {\partial L_2}{\partial \lambda_1}(0, 0, 0) &=  b  c_0\int_0^{\a} e^{-y} \phi_1(y) \, dy, \\
    \frac {\partial L_2}{\partial \lambda_2}(0, 0, 0) &=  b  c_0\int_0^{\a} e^{-y} \phi_2(y) \, dy,
\end{aligned}
\ee
which are all positive by the choice of $\phi_1$ and $\phi_2$. Recall that $\phi_2(x) = \phi_1(x - \a/3)$. Plugging this back into \eqref{b15} and applying a change of variables, we get
\be
\label{b16}
\begin{aligned}
   & \frac {\partial L_1}{\partial \lambda_1}(0, 0, 0)= a_{11}, \,
    \frac {\partial L_1}{\partial \lambda_2}(0, 0, 0)= e^{\a/3} a_{11}, \,\\
   & \frac {\partial L_2}{\partial \lambda_1}(0, 0, 0) = a_{21}, \,
    \frac {\partial L_2}{\partial \lambda_2}(0, 0, 0) = e^{-\a/3} a_{21}. 
    \end{aligned}
\ee

Plugging this into \eqref{b14} then gives us

\be
\label{b17}
\det \begin{bmatrix}
    \frac {\partial L_1}{\partial \lambda_1} & \frac{\partial L_1} {\partial \lambda_2} \\
    \frac {\partial L_2}{\partial \lambda_1} & \frac{\partial L_2} {\partial \lambda_2}
\end{bmatrix}(0,0,0) = (e^{-\a/3}- e^{\a/3})a_{11}a_{21} \ne 0.
\ee

Therefore, there is a neighborhood  $(-\delta,\delta)$, of $0\in \R$ and smooth functions $h_1, h_2:(\delta,\delta)\to \R$ with
$h_1(0)=h_2(0)=0$ such that $L(h_1(\lambda_3),h_2(\lambda_3),\lambda_3)=(0,0)$ for $\lambda_3\in (-\delta,\delta)$. By choosing $\lambda_3 \ne 0$ inside this neighborhood, we find our desired $u_0 = c_0 + h_1(\lambda_3) \phi_1 +h_2( \lambda_3) \phi_2 + \lambda_3 \phi_3$. This completes the proof.

\end{proof}
\vskip.1in



\bigskip

\begin{thebibliography}{99}


\bibitem{BSS} R.  Beals, D. Sattinger, and J. Szmigielski, {\em Multi-peakons and the classical moment problem},  Adv. Math. {\bf 154} (2000),  
229--257.

 \bibitem{BSS2} R.  Beals, D. Sattinger, and J. Szmigielski, {\em Multi-peakons and a theorem of Stieltjes},
Inverse Problems 15 (1999), no. 1, L1--L4. 

\bibitem{Be}  T. B. Benjamin, \emph{Internal waves of permanent
form in fluids of great depth},
J. Fluid Mech. {\bf 29} (1967), 559--592.
{\bf 31} (2013) 2--53.




 \bibitem{BBM} T. B. Benjamin, J. L. Bona, and J. J. Mahony, {\em Model equations for long waves in nonlinear dispersive systems},
Philos. Trans. Roy. Soc. London Ser. A 272 (1972), no. 1220, 47--78. 



 
 \bibitem{BrCo} L. Brandolese and M. F. Cortez, {\em On permanent and breaking waves in hyperelastic rods and rings},
J. Funct. Anal. {\bf 266} (2014), 6954--6987. 



\bibitem{BC} A. Bressan and A. Constantin, {\em Global conservative solutions to the Camassa-Holm equation},  Arch. Rat. Mech. Anal. {\bf 183} (2007), 215--239.

\bibitem{BCZ} A. Bressan, G. Chen,  and Q. Zhang, {\em Uniqueness of conservative solutions to the Camassa-Holm equation
via characteristics}, Discr. Cont. Dyn. Syst. {\bf 35} (2015), 25--42.







\bibitem{CH} R. Camassa and D. Holm, {\em An integrable shallow water equation with peaked solitons}, Phys. Rev. Lett. {\bf 71} (1993), 1661--1664.

\bibitem{CocKa} G. M. Coclite and K. H. Karlsen {\em On the well-posedness of the Degasperis-Procesi equation},  J. Funct. Anal. {\bf 233} (2006), 60--91.


\bibitem{CoEs1} A. Constantin and J.  Escher, {\em Global existence and blow-up for a shallow water equation},  Ann. Scuola Norm. Sup. Pisa Cl. Sci. {\bf 26} (1998), no. 2, 30--328.

\bibitem{CoEs2} A. Constantin and J. Escher,{\em  Wave breaking for nonlinear nonlocal shallow water equations}, Acta Math. {\bf 181} (1998), no. 2, 229--243.

\bibitem{CoEs3} A. Constantin and J. Escher,{\em  Well-posedness, global existence, and blowup phenomena for  a periodic quasi-linear hyperbolic equation}, Comm. Pure Appl. Math. {\bf 51} (1998), 475--504.


\bibitem{CoMc} A. Constantin and H.  McKean, 
{\em A shallow water equation on the circle},
Comm. Pure Appl. Math. {\bf 52} (1999), 949--982. 


\bibitem{CoMo} A. Constantin and L. Molinet, {\em Global weak solutions for a shallow water equation}, Comm. Math. Phys. {\bf 211} (2000), 45--61.









\bibitem{Dai} H.-H. Dai, {\em Model equations for nonlinear dispersive waves in a compressible Mooney-Rivlin rod}, Acta Mech. {\bf 127} (1998), 193--207.
 
  \bibitem{DaHu} H.-H. Dai and Y. Huo {\em Solitary shock waves and other travelling waves in a general compressible hyper-elastic rod},  Proc. R. Soc. Lond. Ser. A Math. Phys. Eng.  {\bf 456} (2000), 331-363.
  
  \bibitem{DP} A. Degasperis, M. Procesi, {\em Asymptotic integrability}, Symmetry and Perturbation Theory, World Scientific, Singapore, (1999) 23--37.


\bibitem{LKT} C. de Lellis, T. Kappeler, and P. Topalov, {\em Low-regularity solutions of the periodic Camassa-Holm equation}, 
Comm. PDE  {\bf 32} (2007), no. 1-3, 8--126. 




\bibitem{EY} J. Escher and Z. Yin, {\em Well-posedness, blow-up phenomena, and global solutions for the b-equation}, J. Reine Angew. Math.
 {\bf 624} (2008), 51--80.


  
 \bibitem{FF}  B. Fuchssteiner and A. S. Fokas, {\em  Symplectic structures, their B\"acklund transformations and hereditary symmetries},  Phys. D 4 (1981/82),  4--66.
 

\bibitem{GHR1} K. Grunert, H. Holden, and X. Raynaud, {\em Lipschitz metric for the Camassa-Holm equation on the line},  Discrete Contin. Dyn. Syst. {\bf 33} (2013), no. 7, 2809--2827.
 
\bibitem{GHR2} K. Grunert, H. Holden, and X. Raynaud, {\em Global conservative solutions to the Camassa-Holm equation for initial data with nonvanishing asymptotics},  Discrete Contin. Dyn. Syst. {\bf 32}  (2012), no. 12, 4209--4227.

\bibitem{HHG} A.A. Himonas, C. Holliman, and K. Grayshan, {\em Norm inflation and ill-posedness for the Degasperis-Procesi equation}, Comm. Partial Differential Equations {\bf 39}
(2014), 2198--2215.

\bibitem{HKM} A.A. Himonas, C. E. Kenig, and G. Misio\l ek,  {\em Non-uniform dependence for the periodic CH equation}, Comm. Partial Differential Equations {\bf 35} (2010), 1145--1162.





\bibitem{HoSt} D. D. Holm and M. F. Staley, {\em Nonlinear balance and exchange of stability of dynamics of solitons, peakons,
ramps/cliffs and leftons in a 1 + 1 nonlinear evolutionary PDE} Phys. Lett. {\bf A 308} (203) 437?44

 \bibitem{HP} C. Hong and G. Ponce, \emph{On special properties of solutions to the Benjamin-Bona-Mahony equation},  J. Diff. Eqs., 393, (2024), 321--342.
 
 \bibitem{Iv} R. Ivanov, {\em Water waves and integrability}, Philos. Trans. Roy. Soc. Lond. A {\bf 365} (2007) 2267--2280.
 
  \bibitem{KPV19} C. E. Kenig, G. Ponce, and L. Vega, {\em Uniqueness  properties of solutions to the Benjamin-Ono equation and related models}, J. Funct. Anal. {\bf 278} (2020), 14 pages.
  
\bibitem{KdV} D. J. Korteweg and G. de Vries,
\emph{On the change of form of long waves advancing in a rectangular canal, and on a new type of long stationary waves}, 
Philos. Mag. (5) {39} (1895), no. 240, 422--443.

 \bibitem{Le} J. Lenells, \emph{Traveling wave solutions of the Camassa-Holm equation},  J. Diff. Eqs., 215, (2005), 393--430.

\bibitem{LiOl} Y. A. Li and P. J. Olver, \emph{Well-posedness and Blow-up Solutions for an Integrable Nonlinearly Dispersive Model Wave Equation}, J. Diff. Eqs. {\bf 162} (2002), 27--63.


\bibitem{LiPo} F. Linares and G. Ponce, \emph{Unique continuation properties for solutions to the Camassa-Holm equation and related models},
Proc. Amer. Math. Soc. {148}  (2020), 3871--3879.

\bibitem{LiPoSi} F. Linares, G. Ponce, and T. Sideris, {\em Properties of solutions to the Camassa-Holm equation on the line in
 a class containing the peakons}, Advanced Studies in Pure Math.,  Asymptotic Analysis for Nonlinear Dispersive 
 and Wave Eqs, {\bf 81} (2019), 196--245.




\bibitem{Ma}  Y. Matsuno, {\em Multisoliton solutions of the Degasperis-Procesi equation and their peakon limit}, Inverse Problems, 21  (2005), no 5, 1553--1570.

\bibitem{Mc} H.  McKean, {\em Breakdown of the Camassa-Holm equation}, Comm. Pure  Appl. Math. {\bf 57} (2004), 416--418.

\bibitem{M} L. Molinet, {\em On well-posedness results for Camassa-Holm equation on
the line: a survey},  J. Nonlinear Math. Phys. {\bf 11} (2004), 521--533.




\bibitem{On} H. Ono,  \emph{Algebraic solitary waves in stratified fluids}, Journal  Physical Society of Japan, {\bf 39} (4) (1975), 108--1091.

\bibitem{Par} A. Parker, {\em On the Camassa-Holm equation and a direct method of solutions. II. Solitons solutions}, Proc. R. Soc. A. {\bf 461} (2005), 3611--3632.



\bibitem{Ro} G. Rodriguez,  \emph{On the Cauchy problem for the Camassa-Holm equation}, Nonlinear Analysis {\bf 46} (2001), 309--327.
{7} (2010), 289---305.



\bibitem{SaSc} J-C. Saut and B. Scheurer, {\em Unique continuation for some evolution equations}, J. Diff. Eqs. {\bf 66} (1987), 118--139.





\end{thebibliography}
\end{document}